\documentclass[11pt]{amsart}

\usepackage{amsmath,amssymb,amscd,amsfonts,verbatim}

\newtheorem{thm}{Theorem}[section]

\newtheorem{assu-nota}[thm]{Assumption--Notation}
\theoremstyle{remark}

\newcommand{\Q}{\mathbb Q}

\newcommand{\pp}{\mathbb P}

\newcommand{\Ups}{\Upsilon}
\newcommand{\al}{\alpha}

\newcommand{\De}{\Delta}
\newcommand{\Si}{\Sigma}

\newcommand{\fie}{\varphi}
\newcommand{\Up}{\Upsilon}
\newcommand{\OO}{\mathcal{O}}

\numberwithin{equation}{section}

\begin{document}

\title{The degree of the bicanonical map of a surface with $p_g=0$.}
\author{Margarida Mendes Lopes  and Rita Pardini}
\date{}
\begin{abstract}
In this note it is shown that, given  a  smooth minimal complex surface of general type $S$  with $p_g(S)=0$, $K^2_S=3$,   for which  the bicanonical map  $\fie_{2K}$ is a morphism,  the degree of $\fie_{2K}$ is not 3. This completes our earlier results, showing that if $S$ is a minimal surface of general type with $p_g=0$, $K^2\ge 3$ such that $|2K_S|$ is free, then the bicanonical map  of $S$ can have degree 1, 2 or 4.

\end{abstract}
\maketitle
\section{Introduction}
Complex surfaces of general type with $p_g=0$ continue to be intriguing.  Continuing our study of their bicanonical map,  
 in this  note  we prove the following:
 
\begin{thm}\label{degree3}Let $S$ be a smooth minimal complex surface of general type with $p_g(S)=0$ and $K^2_S=3$ such that the bicanonical system $|2K_S|$ is base point free. Then the degree of the bicanonical map of $S$ is different from 3.
\end{thm}

Theorem \ref{degree3}, together with previous results (\cite{marg}, \cite{london1}, \cite{K6})  gives the following general statement on the degree of the bicanonical map of a surface of general type with $p_g=0$ and base point free bicanonical system.

\begin{thm}\label{general} Let $S$ be a smooth minimal complex surface of general type with $p_g(S)=0$, let $\fie\colon S\to\pp^{K^2}$ be the bicanonical map. Then:
\begin{enumerate}
\item if $K^2_S=9$, then $\fie$ is birational;
\item if $K^2_S=7,8$, then $\deg\fie \le 2$;
\item if $3\le K^2_S\le 6$ and $|2K_S|$ is base point free, then $\deg\fie$ is equal to either $1,2$ or $4$.
\end{enumerate}
\end{thm}

The proof of Theorem \ref{degree3}  is based on  some properties of  non normal rational  quartic surfaces of $\pp^3$, which are studied in the next section.

\medskip
\noindent{\bf Notation and conventions.} We work over the complex numbers; all varieties are
assumed to be compact and algebraic.

 We do not distinguish between line bundles and divisors on a smooth variety. Linear equivalence of divisors is
denoted by
$\equiv$ and numerical equivalence by $\sim$. The remaining notation is standard in algebraic
geometry.

\medskip

\noindent{\bf Acknowledgements.} 
 The first 
author is a member of the Center for Mathematical Analysis,
Geometry and Dynamical Systems, IST  and was  partially supported
by FCT (Portugal) through program  POCTI/FEDER and Project
POCTI/MAT/44068/2002.
The second author is a member of  G.N.S.A.G.A. and was partially supported by P.R.I.N 2002 and 2004 ``Geometria sulle variet\`a' algebriche'' of M.I.U.R.

\section{Non normal rational quartics of $\pp^3$.}\label{quartics}

In this section we show the following:

\begin{thm}\label{fibration} 
Let $\Si\subset \pp^3$ be a  rational quartic surface such that $\Si$ is not normal, i.e. such that the singular locus of $\Si$ has dimension 1. Let $\Upsilon\to\Si$ be the minimal desingularization and  let $|H|$ be the pull back on $\Up$ of the linear system of planes of $\pp^3$. If the  linear system  $|H|$ is complete, then: 

\begin{enumerate}

\item the genus $g(H)$ of $H$ is equal to $2$;

\item  $K_{\Up}^2=0$ and the linear system $|K_{\Up}+H|$ is a free pencil  of rational curves.

\end{enumerate}
\end{thm}

\begin{proof}
We start by noticing that,  since the linear system $|H|$ is the pull back  on $\Upsilon$ of the linear system of planes of $\pp^3$ and $\Si$ has degree 4, $H^2=4$ . The general curve $H$  of $|H|$ is smooth and irreducible by Bertini's theorem
and it is mapped birationally onto a singular plane quartic, thus we have $g(H)=1+(K_{\Up}+H)H/2\le 2$. On the other hand, by the regularity of $\Ups$ the restriction of $|H|$ to a general $H$ is a complete system, hence Riemann--Roch on $H$ gives $g(H)\ge 2$. This proves  (i).

We divide the proof of  assertion (ii) into steps. 
Consider the linear system $|D|\!:=|K_{\Up}+H|$. By (i), we have $K_{\Up}H=-2$ and thus  $DH=2$, $D^2=K_{\Up}^2$. Using the adjunction sequence for a general $H$ one sees that $h^0(D)=2$, namely $|D|$ is a pencil.

Write $|D|=F+|M|$, where $F$ is the fixed part of $|D|$ and $|M|$ is the moving part. 
\medskip

\noindent{\bf Step 1:} {\em $D$ is nef. In particular, we have $D^2\ge 0$.}\newline
 Since $q(\Up)=0$,  the restriction of $|D|$ to the general curve of $|H|$ is the complete canonical system. So for any irreducible component $\theta$ of $F$ we have $\theta H=0$ and, by the index theorem, $\theta^2<0$. 
 
 Let $\theta$ be an irreducible  curve such that $\theta D<0$. Since $D$ is effective, $\theta$ is a component of $F$. Hence $\theta H=0$, $\theta K_{\Up}<0$, $\theta^2<0$, namely $\theta$ is a $-1-$curve contracted by $|H|$, against the assumption that $\Up\to \Si$ is the minimal desingularization.
 \medskip

 \noindent{\bf Step 2:} {\em One has $D^2\le 1$. If $D^2=1$, then $H=-2K_{\Up}$.}\newline
  The inequality $D^2\le 1$ follows by the index theorem applied to $D$ and $H$. If $D^2=1$, the index theorem gives also $H\sim 2D$, i. e. $H\sim -2K_{\Up}$. Since $\Up$ is rational, we actually have $H\equiv -2K_{\Up}$.
   \medskip
   
   \noindent{\bf Step 3:} {\em The case $D^2=1$ does not occur.}\newline
   By Step 2, if $D^2=1$, then $-K_{\Up}$ is nef and $K^2_{\Up}=1$. Riemann--Roch gives $h^0(-K_{\Up})=2$. Write $|-K_{\Up}|=\Delta+|N|$, where $\De$ is the fixed part and $|N|$ is the moving part. We have:
   $$1=K^2_{\Up}=-K_{\Up}\Delta -K_{\Up}N\ge -K_{\Up}N= N^2+\De N\ge 0.$$
  Assume that $\De\ne 0$. Then, since a nef and big divisor is 1-connected, we have $\De N>0$. Hence in this case all the previous inequalities are equalities and we have
    $$\De N=1, \ N^2=0,\ \  K_{\Up}\De=0.$$
It follows $\De^2=-K_{\Up}\De-N\De=-1$, contradicting the fact  that $\De^2+K_{\Up}\De$ is even  by the adjunction formula.  So we conclude that  $\De=0$.  Since $K_{\Up}^2=1$, it follows that the general curve of $|-K_{\Up}|$ is smooth and irreducible of genus 1. Since $(-K_{\Up})H=2$, the restriction of $|H|$ to the general curve of $|-K_{\Up}|$ is not birational, contradicting the assumption that $|H|$ is the pull back of a very ample system via a birational morphism. Hence the case $D^2=1$ does not occur.
 \medskip

 \noindent{\bf Step 3:} {\em $D^2=0$ and  $|D|$ is a pencil of rational curves.}\newline
We have $D^2=K_{\Up}^2=0$ by the previous steps.   Since $D$ and $M$ are nef, we have $$0=D^2\ge DM=M^2+MF\ge MF\ge 0.$$
It follows $M^2=MF=F^2=0$. Hence, by Zariski's lemma (see \cite{BPV}, Lemma (8.2), p. 90) there is a rational number  $\al\ge 0$ such that $F\equiv \al M$ as $\Q-$divisors. One has:  $2=-K_{\Up}D=-K_{\Up}M(1+\al)$.
Since $-K_{\Up}M=\frac{2}{1+\alpha}>0$,  by the adjunction formula we must have $-K_{\Up}M=2$, $\al=0$.
Hence $F=0$ and $|D|$ is a free pencil of rational curves.
\end{proof}

\section{The proofs of Theorem \ref{general} and Theorem \ref{degree3}}

\begin{proof}[Proof of Theorem \ref{degree3}]
Since $K_S^2=3$, $h^0(S, 2K_S)=K_S^2+\chi(\OO_S)=4$. 
Let $\fie\colon S\to\pp^3$ be the bicanonical map and 
denote by $\Si\subset \pp^3$ the image of $\fie$, which is a surface by \cite{xiao}.  

We argue by contradiction. So from now on we assume that $\deg\fie=3$.
Since $|2K_S|$ is free, we have:
$$4K^2_S=12=\deg\fie\deg\Si=3\deg\Si,$$
namely $\Si$ is a quartic surface.

Let $\Up\to \Si$ be the minimal desingularization. In principle, $\Up$ is either a K3 surface or a ruled surface.  Since $\Up$ is birationally dominated by $S$ and $q(S)=0$ by the inequality $\chi(\OO_S)>0$ for surfaces of general type, we have  $p_g(\Up)=q(\Up)=0$ and so $\Up$ is rational. We note also that, since $\Si$ is the image of $S$ by a complete linear system, the linear system $|H|$, as in the assumptions of Theorem \ref{fibration}, is   complete.

We wish to show that $\Si$ is not normal.
Assume by contradiction that $\Si$ is normal. Being a quartic hypersurface, $\Si$ is a  Gorenstein surface with $K_{\Si}=0$ and $p_g(\Si)=1$, namely there exists a nonzero 2-form $w$  on $\Si$ which is regular outside the singular locus of $\Si$, which is a finite set. The pullback $w'$ of $w$ to  $S$ is a nonzero 2-form which is regular on the complement of the union of the $-2-$curves of $S$ and of a finite set of points. It is well known (cf. \cite{steenbrink}, Lemmas 1.8 and 1.11) that such a form $w'$ is actually regular on $S$, but this contradicts the assumption $p_g(S)=0$. So $\Si$ is not normal and  by  Theorem \ref{fibration}  we have $g(H)=2$ and $K_{\Up}H=-2$.

Consider the linear system  $|K_{\Up}+H|$. By Theorem \ref{fibration}, this system defines a fibration $f\colon \Up\to \pp^1$ whose general fibre  is a rational curve. Since $K_{\Up}^2=0$ by Theorem \ref{fibration}, the fibration $f$ is not relatively minimal. Let $\theta$ be an irreducible $-1-$curve contained in a fibre of $f$. We have $\theta(K_{\Up}+H)=0$, namely $\theta H=1$. 
So the image of $\theta$ in $\Si$ is a line $L$. If $L'$ is the  pull back  of $L$ on $S$, then we have
$3=L'\psi^*H=L'(2K_S)$, a contradiction.

\end{proof}

\begin{proof}[Proof of Theorem \ref{general}] 
Statement (i) and (ii) are the main result of \cite{london1}.

So we consider the case $3\le K^2_S\le 6$.
 We have $\deg\fie\le 4$, by \cite{marg}. Since $|2K_S|$ is free by assumption, we have:
$$4K^2_S=\deg\fie\deg\Si.$$
Hence we have to exclude the cases $K^2_S=6$, $\deg\fie=3$ and $K^2_S=\deg\fie=3$. The former does not occur by \cite[Thm.1.1]{K6} and the latter is excluded by Theorem \ref{degree3}.

\end{proof}

\bigskip

\bigskip

\begin{minipage}{13cm}
\parbox[t]{6.5cm}{Margarida Mendes Lopes\\
Departamento de  Matem\'atica\\
Instituto Superior T\'ecnico\\
Universidade T{\'e}cnica de Lisboa\\
Av.~Rovisco Pais\\
1049-001 Lisboa, Portugal\\
mmlopes@math.ist.utl.pt
 } \hfill
\parbox[t]{5.5cm}{Rita Pardini\\
Dipartimento di Matematica\\
Universit\`a di Pisa\\
Largo B. Pontecorvo, 5\\
56127 Pisa, Italy\\
pardini@dm.unipi.it}
\end{minipage}

\end{document}